\documentclass[12pt]{amsart}

\newtheorem{theorem}{Theorem}[section]
\newtheorem{defn}{Definition}[section]
\newtheorem{lem}{Lemma}[section]

\allowdisplaybreaks[4]
\begin{document}

\title{Extensions of Weak-Type Multipliers}


\author{P. Mohanty}


\address{Department of Mathematics, Indian Institute of Technology, 
Kanpur-208016, India}




\email{parasar@iitk.ac.in}


\thanks{The first author was supported by CSIR}


\author{S. Madan}

\address{Department of Mathematics, Indian Institute of Technology, 
Kanpur-208016, India}

\email{madan@iitk.ac.in}







\begin{abstract}
 In this paper we prove that if $\Lambda\in  M_p(\mathbb R^N)$ and has 
 compact support then $\Lambda$ is a weak summability kernel for 
 $1<p<\infty$, where $  M_p(\mathbb R^N)$ is the space of multipliers 
 of $L^p(\mathbb R^N)$.
\end{abstract}


\maketitle

\clearpage

\section{Introduction}

 Let $G$ be a locally compact abelian group, with Haar measure 
 $\mu$ and
 let $\hat G$ be its dual. We call an operator 
$T: L^p(G)\longrightarrow L^{p,\infty}(G),\;1\leq p<\infty$, a multiplier of weak type 
$(p,p)$, 
 if it is bounded and translation invariant i.e.
$\tau_x T=T\tau_x \;\;\;\forall x\in G$, and there exists a constant 
$C>0$ such that 
\begin{eqnarray}
\mu\{x\in G : |Tf(x)|>t\}\leq \frac{C^p}{t^p}\|f\|_p^p 
\label{1.1}
\end{eqnarray}
for all 
$f\in L^p(G)$ and $t>0$. (Here $L^{p,\infty}$ denotes the standard weak $L^p$ spaces.) 
 Asmar, Berkson and Gillespie 
in \cite{ABG3} proved that for 
all such operators $T$ there exists a $\phi\in L^\infty(\hat G)$ 
such that 
$(Tf)^\wedge =\phi \hat f$ for all $f\in L^2\cap L^p(G)$. We will also 
call such $\phi$'s to be multipliers of weak type $(p,p)$. Let  
$M_p^{(w)}(\hat G)$ denote the space of multipliers of 
weak type $(p,p)$ for $1\leq p<\infty$, and let $  N_p^{(w)}(\phi)$ be 
the  smallest constant $C$ such that inequality (\ref{1.1}) holds.

In this paper we are concerned with extensions of weak type multipliers 
from $\mathbb Z^N$ to $\mathbb R^N$ through summability kernels. For 
similar results on strong type 
multipliers , see \cite{J}, \cite{BPW}. 
Here we identify $\mathbb T^N$ with  $[0, 1)^N$ 
and for $f\in {L^1(\mathbb R^N)}$ we define its Fourier transform as
$\hat{f}(\xi)=\int\limits_{\mathbb R^N}f(x)\;e^{-2\pi i\xi\cdot x}dx$ 
for $\xi\in\mathbb R^N$. Let us 
define summability kernels for weak type multipliers as follows
\begin{defn}
 A bounded measurable function 
 $\Lambda :\;\mathbb R^N\longrightarrow \;\mathbb C$ is 
 called a weak summability kernel for $  M_p^{(w)}(\mathbb R^N)$ if 
 for $\phi\in  M_p^{(w)}(\mathbb Z^N)$ the function 
$W_{\phi, \Lambda}(\xi)= \sum\limits_{n\in\mathbb Z^N}\phi(n)
\Lambda(\xi-n)$  is defined and belongs
to $  M_p^{(w)}(\mathbb R^N)$. 
\end{defn}
This definition is just the weak type analouge of summability kernel 
for 
strong type multipliers~\cite{BPW}. We first cite two important results regarding
 the summability kernels of strong type multipliers from the work of 
Jodeit \cite{J} and of Berkson, Paluszynski and Weiss \cite{BPW}:
\begin{theorem}{\rm\cite{J}}
Let $S\in L^1(\mathbb R^N)$ and 
$supp\;S\subseteq[\frac{1}{4},\frac{3}{4}]^N$ with
$\tau= \sum\limits_{n\in\mathbb Z^N} |\hat s(n)|<\infty$ , where 
$s$ is the 1-periodic extension of 
 $S$, then the function defined by 
$W_{\phi,\hat S}(\xi)= \sum\limits_{n\in\mathbb Z^N}\phi(n)\hat 
S(\xi-n)$
 belongs to $M_p(\mathbb R^N)$, for $1\leq p<\infty$ with
 $\|W_{\phi,\hat S}\|_{M_p(\mathbb R^N)}\leq 
C_p\tau\|\phi\|_{M_p(\mathbb Z^N)}$
\end{theorem}
\begin{theorem}{\rm \cite{BPW}}
For $1\leq p<\infty$, let $\Lambda\in M_p(\mathbb R^N)$ and 
$supp \Lambda\subseteq [\frac{1}{4},\frac{3}{4}]^N$. For 
$\phi\in M_p(\mathbb Z^N)$ define 
$W_{\phi, \Lambda}(\xi)= \sum\limits_{n\in\mathbb Z^N}\phi(n)
\Lambda(\xi-n)$ on $\mathbb R^N$. Then 
$W_{\phi, \Lambda}\in  M_p^(\mathbb R^N)$ and
 $\|W_{\phi, \Lambda}\|_{M_p(\mathbb R^N)}\leq C_p\|
 \Lambda\|_{M_p(\mathbb R^N)}\|\phi\|_{M_p(\mathbb Z^N)}$
 where  $C_p$ is a constant. (Further, if $\Lambda$ has arbitary compact 

support the same result holds except that the constant $C_p$ necessarily depends 

on the support of $\Lambda$, as shown in ~\cite{BPW}
\label{thm:9.4}
\end{theorem}
Asmar, Berkson and Gillespie proved a weak type analogue of 
Theorem 1.1 in 
\cite{ABG6}.  In this same paper they also proved that 
$\Lambda$ defined by  
$\Lambda(\xi)=\prod\limits_{j=1}^N \max(1-|\xi_j|,0)$ for 
$\xi=(\xi_1,...,\xi_N)$ is a weak type summability kernel. 
 In this paper, we prove the weak type analouge of Theorem 1.2 in 
 $\S 2$, for $1<p<\infty$. In $\S 3$ we relax the hypothesis that 
 $supp~\Lambda\subseteq [\frac{1}{4},\frac{3}{4}]^N$. For the proof of 
our main result , as in \cite{BPW}, we will obtain the weak type 
inequalities 
by applying the technique of transference couples due to Berkson, 
Paluszy${\rm{\check n}}$ski,  and Weiss~\cite{BPW}. 
\begin{defn}
 For a locally compact group $G$, a transference couple is a pair 
$(S,T)=(\{S_u\},\{T_u\})$ , $u\in G$, of 
strongly continuous mappings defined on $G$ with values in 
${\mathcal B}(X)$, 
where $X$ is a Banach space, satisfying\\
(i) $C_S=\sup\{\|S_u\|:u\in G\}<\infty$\\
(ii) $C_T=\sup\{\|T_u\|:u\in G\}<\infty$\\
(iii) $S_v T_u=T_{vu}\;\;\;\forall u,v\in G$
\end{defn}
  In $\S 4$, as an application of our result, we prove a weak type 
analogue of an extension theorem by 
de Leeuw.
\section{Weak-Type Inequality for Transference Couples and The 
Main Theorem}
Let $\Lambda\in L^\infty(\mathbb R^N)$ and $supp\;\Lambda\subseteq 
[\frac{1}{4},\frac{3}{4}]^N$. 
Consider the following  transference couple $(S,T)$ used by 
Berkson, Paluszy$\rm{\acute n}$ski, and Weiss in~\cite{BPW}. For 
$u\in\mathbb T^N$ the family 
$T=\{T_u\}$ is given by 
\begin{equation}
(T_u f)^\wedge(\xi)= \sum\limits_{n\in\mathbb Z^N}\Lambda(\xi-n)e^{2\pi 
iu.n} \hat f(\xi),\; \;{\rm for}\;f\in L^p(\mathbb R^N)
\label{eqn:9.4}
\end{equation}
and the family $S=\{S_u\}$ is defined by
\begin{equation}
(S_u f)^\wedge(\xi)= \sum\limits_{n\in\mathbb Z^N} b(\xi-n)e^{2\pi 
iu.n} \hat f(\xi).\; \;{\rm for}\;f\in L^p(\mathbb R^N),
\label{eqn:10.4}
\end{equation}

 where $b(\xi)=\prod\limits _{i=1}^Nb_i(\xi_i)$ for 
 $\xi=(\xi_1,......,\xi_N)$ and
for each $i$, $b_i$ is the continuous function defined on 
$\mathbb R$ as 
$b_i(x)=1$ if $x\in[\frac{1}{4},\frac{3}{4}]$, linear in 
$[0,\frac{1}{4})\cup(\frac{3}{4},1]$ and $0$ otherwise.
 It is easy to see that
\begin{eqnarray}
S_u f(x)= \sum\limits_{l\in\mathbb Z^N}\check \beta 
_u(l)f(x+u-l)\;\;\;a.e.,
\label{eqn:11.4}
\end{eqnarray}
where $\check \beta_u$ is the inverse Fourier transform of the function
 $\beta_u(\xi)=b(\xi)e^{2\pi i\xi.u}$, given explicitily by
 $$\check\beta_u(\xi)=\prod_{i=1}^N\check\beta_{u_i}(\xi_i),$$
where 
\begin{equation}
\check\beta_{u_i}(\xi_i)=
\begin{cases}
  \frac{2e^{2\pi i\frac{(\xi_i+u_i)}{2}}}{\pi^2 
(\xi_i-u_i)^2}(cos\frac{\pi}{2}(\xi_i-u_i)-cos\pi(\xi_i-u_i))
\; & {\mbox{if}}\;\; \xi_i\not= u_i\\
\frac{3e^{2\pi i\frac{(\xi_i+u_i)}{2}}}{4}\; \; & {\mbox{if}} 
\;\xi_i=u_i.
\end{cases}
\label{betau}
\end{equation}
Then by a straightforward calculation using Eqn.(\ref{betau}) we have
\begin{equation}
  \sum\limits_{l\in\mathbb Z^N}|\check\beta_u(l)|\leq  
\sum\limits_{l\in\mathbb Z^N}\beta(l)=C<\infty,
\label{eqn:12.4}
\end{equation}
where $\beta(l)=\prod_{i=1}^N\beta_i(l_i)~{\rm and}$
\begin{eqnarray*}
\beta_i(l_i) = 
\begin{cases}
\frac{1}{(l_i-1)^2} &if\;l_i>1\\
\frac{1}{(l_i+1)^2} &if\;l_i<1 \\ 
\|b_i\|_1 &{\mbox{otherwise.}}
\end{cases}
\end{eqnarray*}
 In the following theorem we shall show that the  operator transferred 
by $T$ (of the transfernce couple $(S,T)$ defined in 
Eqn.~(\ref{eqn:9.4}) and Eqn.~(\ref{eqn:10.4})) given  by
$$H_k f(.)=\int_{\mathbb T^N}k(u)T_{u^{-1}}f(.)du,$$
where  $k\in L^1(\mathbb T^N)$ and $f\in L^p(\mathbb R^N)$, satisfies a
 weak  $(p,p)$ inequality.
\begin{theorem}
 Let $(S,T)$ be the transference couple as defined in 
 Eqn.~(\ref{eqn:9.4}) and Eqn.~(\ref{eqn:10.4}). Then
 for $1<p<\infty$ and $t>0$\\
\[\lambda\{x\in\mathbb R^N\;:\;|H_k f(x)|>t\}\leq 
(\frac{C\;C_p}{t}C_T  N_p^{(w)}(k)\|f\|_p)^p,\]
where $\lambda$ denotes the Lebesgue measure of $\mathbb R^N$, 
$C= \sum\limits_{l\in\mathbb Z^N} \beta(l)$ as in 
Eqn.~(\ref{eqn:12.4}), $C_T$ is the uniform bound for the 
family $T=\{T_u\}$, and $C_p=\frac{p}{p-1}$.
\label{thm:10.4}
\end{theorem}
\noindent{\bf Proof:}
Assume $f\in {\mathcal S}(\mathbb R^N)$. For $t>0$ define 
$E_t=\{x\;:\;|H_k f(x)|>t\}$.

 \noindent Notice that 

$H_k f(x)= S_{v^{-1}}S_v H_k f(x)=\sum\limits_{l\in\mathbb Z^N}\check\beta_{v^{-1}}(l)\int_
{\mathbb T^N}k(u)T_{u^{-1}v}f(x-v-l) du|>t\}.$
Let ${\mathcal F}_t=\{(v,x)\in{\mathbb T^N\times\mathbb R^N}\;:\;| 
\sum\limits_{l\in\mathbb Z^N}\check\beta_{v^{-1}}(l)\int_
{\mathbb T^N}k(u)T_{u^{-1}v}f(x-l) du|>t\}.$

Then, using translation invariant of Lebesgue measure
\begin{eqnarray*}
\lambda(E_t)& = &\lambda\{x\in\mathbb R^N\;:\;|S_{v^{-1}}\int_{\mathbb 
T^N}k(u)T_{u^{-1}v} f(x)du|>t\}\\
& = &\lambda\{x\in\mathbb R^N:| \sum\limits_{l\in\mathbb 
Z^N}\check\beta_{v^{-1}}(l)\int_{\mathbb T^N}k(u)T_{u^{-1}v} f(x-l) du|>t\}\\
& = &\int_{\mathbb T^N}\int_{\mathbb R^N}\chi_{{\mathcal F}_t}(v,x) dx 
dv\\
& = &\int_{\mathbb R^N}|\{v:| \sum\limits_{l\in\mathbb 
Z^N}\check\beta_{v^{-1}}(l)\int_{\mathbb T^N}k(u)T_{u^{-1}v} f(x-l) du|>t\}|dx,
\end{eqnarray*}
 where $|E|$ denotes the measure of the subset $E\subseteq\mathbb T^N$. 
 Thus
\begin{eqnarray*}
\lambda(E_t)& \leq & \int_{\mathbb R^N}|\{v: \sum\limits_
{l\in\mathbb Z^N}\beta(l)|\int_{\mathbb  T^N}k(u)T_{u^{-1}v} 
f(x-l)| du>t\}|dx\\
& = &\int_{\mathbb R^N}|\{v: \sum\limits_{l\in\mathbb 
Z^N}\beta(l)|k*F(.,x-l)(v)|>t\}|dx ,\;\;\; {\rm where}\;F(v,x)=T_v f(x)~~~{\rm a.e.}.
\end{eqnarray*}
 We know that $\sup\limits_{t>0}~ 
t{\lambda_f(t)}^{\frac{1}{p}}=\|f\|_{L^{p,\infty}}$
 for $f\in L^{p,\infty}$. Also, since $p>1$, $\|\;\|_{p,\infty}$ is 
equivalent to a norm $\|\;\|_{p,\infty}^*$ (\cite{SW}), using traingle inequality for norms 

we have 
\begin{eqnarray*}
\lambda(E_t)& \leq &\int_{\mathbb R^N}\frac{1}{t^p}\| 
\sum\limits_{l\in\mathbb Z^N} \beta(l)|k*F(.,x-l)\|_{L^{p,\infty}(\mathbb T^N)}^p dx\\
& \leq &C_p^p\int_{\mathbb R^N}\frac{1}{t^p}( \sum\limits_{l\in\mathbb 
Z^N} \beta(l)\|k*F(.,x-l)\|_{L^{p,\infty}(\mathbb T^N)}^*)^p dx, 
\;\;\;{\mbox{where}}\;\;C_p=\frac{p}{p-1}\\
& \leq &C_p\int_{\mathbb R^N}\frac{1}{t^p}( \sum\limits_{l\in\mathbb 
Z^N} \beta(l)  N_p^{(w)}(k)\|F(.,x-l)\|_{L^p(\mathbb T^N)})^p dx,
\end{eqnarray*}
where $  N_p^{(w)}(k)$ is the weak-type norm of the convolution 
operator 
 $f\longmapsto k*f$ for $f\in L^p(\mathbb T^N)$. Thus,
\begin{eqnarray*}
\lambda(E_t)& \leq &C_p^p\frac{1}{t^p} \sum\limits_{l\in\mathbb 
Z^N}\beta(l)  N_p^{(w)}(k)(\int_{\mathbb R^N}\int_{\mathbb T^N}|T_v f(x-l)|^p 
dx dv)^{\frac{1}{p}})^p\\
& = &C_p^p\frac{1}{t^p}( \sum\limits_{l\in\mathbb Z^N}\beta(l)  
N_p^{(w)}(k)(\int_{\mathbb T^N}\int_{\mathbb R^N}|T_v f(x-l)|^p dx 
dv)^{\frac{1}{p}})^p\\
& \leq &(\frac{C C_pC_T}{t^p}  N_p^{(w)}(k)\|f\|_p)^p.
\end{eqnarray*}
Hence, $H_k f$ satisfies a weak $(p,p)$ inequality.

 In order to prove the weak-type analogue of Theorem~\ref{thm:9.4} we 
 need the 
following Lemma proved by Asmar, Berkson, and Gillespie in \cite{ABG5}. 
\begin{lem}{\rm{\cite{ABG5}}}
 Suppose that $1\leq p<\infty$ , $\{\phi_j\}\subseteq M_p^{(w)}(\hat 
G)$; 
$\sup\{|\phi_j(\gamma)|:j\in \mathbb N,\gamma\in\hat G\}~~<\infty$ and 
 suppose $\phi_j$ converges pointwise a.e. on $\hat G$ to a function 
 $\phi$ . 
If $\liminf\limits_j  N_p^{(w)}(\phi_j)<\infty$ then $\phi\in 
M_p^{(w)}(\hat G)$ 
and $  N_p^{(w)}(\phi)\leq \liminf\limits_j  N_p^{(w)}(\phi_j)$. 
\label{lem:2.4}
\end{lem}

In the following theorem, we use the family of operators $\{T_u\}$ 

defined in (~\ref{eqn:9.4}) with $\Lambda\in M_p(\mathbb R^N$) and 

$supp\Lambda\subseteq [\frac{1}{4},\frac{3}{4}]^N$. In ths case, 

by ~\cite{BPW} we have $C_T\leq c_p\|\Lambda\|_{M_p(\mathbb R^N})$, 

where $c_p$ is a constant.
\begin{theorem}
 Suppose $1<p<\infty$ and $\Lambda\in   M_p(\mathbb R^N)$  is 
supported
  in the set
 $[\frac{1}{4},\frac{3}{4}]^N$. For $\phi\in  M_p^{(w)}(\mathbb Z^N)$ 
 define
$$W_{\phi, \Lambda}(\xi)= \sum\limits_{n\in\mathbb 
Z^N}\phi(n)\Lambda(\xi-n) \;\;\;\;{\mbox{on}} \;\mathbb R^N.$$
Then $W_{\phi, \Lambda}\in  M_p^{(w)}(\mathbb R^N)$ and $  
N_p^{(w)}(W_{\phi, \Lambda})\leq C  N_p^{(w)}(\phi)\|\Lambda\|_{M_p(\mathbb R^N)}$.
\label{eqn:11a.4}
\end{theorem}
\noindent{\bf Proof:}
Using Lemma~\ref{lem:2.4} we first show that it is enough to prove  the
  theorem for  $\phi\in  M_p^{(w)}(\mathbb Z^N)$ 
having finite support. Suppose the theorem is true for finitely 
supported $\phi$. 
Then for  arbitrary $\phi\in  M_p^{(w)}(\mathbb Z^N)$, 
define $\phi_j=\hat k_j\phi$, where $k_j$ is the j-th
F$\acute{\mbox{e}}$jer  kernel. Then for each $j$ , $\phi_j$'s have 
finite support and 
$(T_{\phi_j}f)^\wedge(n)=\phi_j(n)\hat f(n)=(T_\phi(k_j*f))^\wedge(n)$. 
So $\phi_j\in  M_p^{(w)}(\mathbb Z^N)\;\;$ for each $j$ and 
$  N_p^{(w)}(\phi_j)\leq   N_p^{(w)}(\phi)$. Define 
$W_{\phi_j, \Lambda}(\xi)= \sum\limits_{n\in\mathbb 
Z^N}\phi_j(n)\Lambda(\xi-n)$. Now $\liminf\limits_jW_{\phi_j, \Lambda}(\xi)=W_{\phi, 
\Lambda}(\xi)$.
 Also, by our assumption 
\begin{eqnarray*}
  N_p^{(w)}(W_{\phi_j, \Lambda})& \leq & C  
N_p^{(w)}(\phi_j)\|\Lambda\|_{M_p(\mathbb R^N)}\\
& \leq &C  N_p^{(w)}(\phi)\|\Lambda\|_{M_p(\mathbb R^N)}
\end{eqnarray*}
and $|W_{\phi_j, \Lambda}|\leq 2\|\Lambda\|_\infty\|\phi_j\|_\infty\leq 
2\|\Lambda\|_\infty\|\phi\|_\infty$.
 Thus by Lemma~\ref{lem:2.4}, applied to $W_{\phi_j,\Lambda}$'s ,we 
 conclude that
 $W_{\phi, \Lambda}\in  M_p^{(w)}(\mathbb R^N)$. Hence it is enough 
 to assume that $\phi\in   M_p^{(w)}(\mathbb Z^N)$ 
has finite support.

Now let $\phi\in M_p^{(w)}(\mathbb Z^N)$ be finitely supported. Define
 $k(u)= \sum\limits_{n\in\mathbb Z^N}\phi(n)e^{-2\pi iu.n}$ then 
 $k\in L^1(\mathbb T^N)$ and 
$\hat k(n)=\phi(n)$. For this particular $k$ and the transference 
couple 
$(S,T)$ defined above. We have
$$(H_k f)^\wedge(\xi)=(T_{W_{\phi,\Lambda}}f)^\wedge(\xi).$$
 Thus $T_{W_{\phi,\Lambda}}f=H_kf$. Hence from Theorem~\ref{thm:10.4} 
 and since $C_T\leq c_p\|\Lambda\|_{M_p(\mathbb R^N}$, we have
$$\lambda\{x\in\mathbb R^N:|T_{W_{\phi,\Lambda}}f(x)|>t\}\leq 
(\frac{C}{t}  N_p^{(w)}(\phi)\|\Lambda\|_{M_p(\mathbb 
R^N)}\|f\|_p)^p.$$
\section{Lattice Preserving Linear Transformations and Multipliers}
We shall now relax the hypothesis that  
$supp\;\Lambda\subseteq[\frac{1}{4},\frac{3}{4}]^N$ to allow 
$\Lambda$ to have
 arbitrary compact support. In fact this can be done by a 
partition of identity argument as in \cite{BPW}. Here we give  a 
different method by proving Lemma~\ref{lem:5.4} below.  Particular 
cases of this lemma occur in \cite{J} and in  \cite{ABG6}. 
Suppose $supp\;\Lambda\subseteq [-M,M]^N$; define 
$\Lambda_M(\xi)=\Lambda_1(4M\xi)$, 
where $\Lambda_1(\xi)=\Lambda(\xi-\frac{1}{2})$.  
 So $supp\;\Lambda_M\subseteq [\frac{1}{4},\frac{3}{4}]^N$. Thus if we 
 define a non-singular 
transformation $A:\mathbb R^N\longrightarrow \mathbb R^N$ such that 
$Ax=4Mx$ then 
$\Lambda_M=\Lambda_1\circ A$. In order to replace the support condition 
we need
 to prove $\Lambda_M\circ A^{-1}$ is a summability kernel. In the work 
 of Jodeit and of Asmar, Berkson and  Gillespie they assume $A$ 
in Lemma~\ref{lem:5.4} to be multiplication by 2. We have combined some 
of the results proved by 
  Gr$\ddot{\mbox{o}}$chenig and Madych ~\cite{GM} in the following 
  lemma which will help us to prove 
Lemma~\ref{lem:5.4}. In the proof of Theorem~\ref{thm:12.4}, we only 
use the case of a diagonal linear transform, but the more general 
results proved below are of some interest in their own right.
\begin{lem}{\rm{\cite{GM}}}
 Let $A:\mathbb R^N\longrightarrow \mathbb R^N$ be a non-singular 
linear 
 transformation   which preserves the lattice $\mathbb Z^N$ (i.e.  
$A(\mathbb Z^N)\subseteq \mathbb Z^N$). Then the following are true.\\
(i) The number of distinct coset representatives of 
$\mathbb Z^N/{A\mathbb Z^N}$ is equal to $q=|\det A|$.\\
(ii) If $Q_0=[0,1)^N$ and $k_1,.....,k_q$ are the distinct  coset 
representatives  of $\mathbb Z^N/{A\mathbb Z^N}$ then the sets 
$A^{-1}(Q_0+k_i)$ are mutually disjoint.\\
(iii) Let $Q=\cup _{i=1}^q A^{-1}(Q_0+k_i)$, then $\lambda(Q)=1$ and 
$\cup_{k\in \mathbb Z^N}(Q+k)\simeq \mathbb R^N.$\\
(iv) $AQ\simeq\cup_{i=1}^q (Q_0+k_i)$.\\

Where $E\simeq F$ if $\lambda (F\bigtriangleup E)=0$.
\label{lem:3.4}
\end{lem}
The above result is essentially contained  in \cite{GM}. 
\begin{lem}
 Let $A$ be as in Lemma~\ref{lem:3.4}. Denote $A^t=B$, where $A^t$ is 
 the transpose of $A$. For $\phi\in l_\infty(\mathbb Z^N)$ define
$$\psi(n)=\phi(Bn)$$
and
\begin{center}
$\eta(n)=
\begin{cases}
\phi(B^{-1}n) &n\in B\mathbb Z^N\\
 0 &otherwise.
\end{cases}$
\end{center}
(i) If $\phi\in  M_p(\mathbb Z^N)$ then $\psi, \eta\in  M_p(\mathbb 
Z^N)$
 with multiplier norms 
not exceeding the multiplier norm of $\phi$.\\
(ii) If $\phi\in  M_p^{(w)}(\mathbb Z^N)$ then 
$\psi, \eta\in  M_p^{(w)}(\mathbb Z^N)$ with weak 
multiplier norms not exceeding the weak multiplier norm of $\phi$.
\label{lem:5.4}
\end{lem}
\noindent{\bf Proof:}
(i) For $f\in L^p(Q_0)$, we let $f$ again denote the periodic extension 
to $\mathbb R^N$. Define $Sf(x)=f(Ax)$ , then $Sf$ is also periodic and 
\begin{eqnarray*}
\int_{Q_0}|Sf(x)|^pdx & = 
&\int_{Q_0}|Sf(x)|^p\sum\limits_j\chi_Q(x-j)dx\\
& = &\sum\limits_j\int_{Q_0+j}|Sf(x)|^p\chi_Q(x)dx\\
& = &\int_Q|Sf(x)|^p dx\\
& = &\frac{1}{|\det A|}\int_{AQ}|f(x)|^p dx\\
& = &\frac{1}{q}\sum\limits_{i=1}^q\int_{Q_0+k_i}|f(x)|^p 
dx\;\;\;\;\;\;\;((iv) \;{\mbox{of Lemma~\ref{lem:3.4}}})\\
& = &\int_{Q_0}|f(x)|^p dx.
\end{eqnarray*}
Thus $S$ is an isometry, i.e., 
$\parallel Sf\parallel _{L^p(Q_0)}=\parallel f\parallel _{L^p(Q_0)}$. 
\noindent
Further,  from the orthogonality relations of the characters (Lemma 1, 
\cite{M}) we have 
\begin{eqnarray*}
(Sf)^\wedge(n)=
\begin{cases}
\hat f(B^{-1}n) &if \;n\in B\mathbb Z^N\\
 0 &otherwise.
\end{cases}
\end{eqnarray*}
For $f\in L^p(Q_0)$ we define an  operator $W$ on $L^p(Q_0)$ given by
 $Wf(x)=\frac{1}{q}\sum\limits_{i=1}^q f(A^{-1}(x+~~k_i)),$ $\;\;$ 
where $k_1,\dots,k_q$ are distinct  cosets representations of 
$\mathbb Z^N/ {A\mathbb Z^N}$. 
Then for a trigonometric polynomial $f$,
$$(Wf)^\wedge(n)  =  \hat f(Bn),$$
and so 
\begin{eqnarray*}
(\int_{Q_0}|Wf(x)|^pdx)^{\frac{1}{p}} & = 
&(\int_{Q_0}|\frac{1}{q}\sum\limits_{i=1}^q f(A^{-1}(x+k_i))|^p dx)^{\frac{1}{p}}\\
& \leq 
&\frac{1}{q}\sum\limits_{i=1}^q(\int_{Q_0}|f(A^{-1}(x+k_i))|^pdx)^{\frac{1}{p}}\\
& = 
&\frac{q^{\frac{1}{p}}}{q}\sum\limits_{i=1}^q(\int_{A^{-1}(Q_0+k_i)}|f(x)|^p dx)^{\frac{1}{p}}.
\end{eqnarray*}
Therefore $\parallel Wf\parallel_{L^p(Q_0)}\leq q^
{\frac{1-p}{p}}\parallel f\parallel _{L^p(Q_0)}$, since 

$\int_{Q_0}|f(x)|^p dx=\int_Q|f(x)|^p dx$ as above.  
It is easy to see that
\begin{equation}
ST_\phi W=T_\eta
\label{eqn:13.4}
\end{equation}
and
\begin{equation}
WT_\phi S=T_\psi
\label{eqn:14.4}
\end{equation}
It follows that, if  $\phi\in  M_p(\mathbb Z^N)$  then 
$\|T_\psi f\|\leq C_p\|\phi\|_{M_p(\mathbb Z^N)}\|f\|_{L^p(Q_0)}.$  
 Also $\|T_\eta f\|_{L^P(Q_0)}\leq C_p\|\phi\|_
 {M_p(\mathbb Z^N)}\|f\|_{L^p(Q_0)}$. Hence 
$\psi,\eta\in M_p(\mathbb Z^N)$.

(ii) For $\phi\in  M_p^{(w)}(\mathbb Z^N)$, we need to calculate the 
distribution 
function of $Sf$ and $Wf$. Denote $E_t=\{x\in Q_0 :|Sf(x)|>t>0\}$. Then
\begin{eqnarray*}
|E_t| & = &\int_{Q_0}\chi_{E_t}(x)dx\\
& = &\int_{Q_0}\chi_{\mathbb R_+}(|f(Ax)|-t)dx\\
& = &\frac{1}{q}\int_{AQ}\chi_{\mathbb R_+}(|f(x)|-t)dx\\
& = &\frac{1}{q}\sum\limits_{i=1}^q\int_{Q_0+k_i}\chi_{\mathbb 
R_+}(|f(x)|-t)dx\\
& = &|\{x:|f(x)>t\}|.
\end{eqnarray*}
\noindent
Therefore, 
\begin{equation}
|\{x\in Q_0 :|Sf(x)>t\}|=|\{x\in Q_0 :|f(x)>t\}|
\label{eqn:15.4}
\end{equation}
Also
\begin{eqnarray*}
|\{x\in Q_0:|Wf(x)|>t\}| & = & |\{x\in Q_0 :| \sum\limits_{i=1}^q 
f(A^{-1}(x+k_i))|>tq\}|\\
& \leq &|\{x\in Q_0 : \sum\limits_{i=1}^q|f(A^{-1}(x+k_i))|>tq\}|\\
& = & \sum\limits_{i=1}^q\int_{Q_0}\chi_{\mathbb 
R_+}(|f(A^{-1}(x+k_i))|-t)dx\\
& = & q \sum\limits_{i=1}^q\int_{A^{-1}(Q_0+k_i)}\chi_{\mathbb 
R_+}(|f(x)|-t)dx.
\end{eqnarray*}
Thus
\begin{equation}
|\{x\in Q_0 :|Wf(x)|>t\}\leq q|\{x\in Q_0 :|f(x)|>t\}|.
\label{eqn:16.4}
\end{equation}
From the relations (\ref{eqn:13.4}) - (\ref{eqn:16.4}), we conclude 
that 
$\psi , \eta\in  M_p^{(w)}(\mathbb Z^N)$ whenever 
$\phi\in  M_p^{(w)}(\mathbb Z^N)$ . Also 
$  N_p^{(w)}(\psi)\leq C  N_p^{(w)}(\phi)$ and
 $  N_p^{(w)}(\eta)\leq C  N_p^{(w)}(\phi)$.

 As an application of this Lemma we get the following result regarding
  weak summability kernels.
\begin{lem}
 Let $A$ be as in Lemma~\ref{lem:3.4}. Suppose $\Lambda$ is  a 
 weak (strong) summabiliy kernel then
 $\Lambda\circ B$ and $\Lambda\circ B^{-1}$ are also weak (strong) 
 summability kernels.
\label{lem:6.4}
\end{lem}
\noindent{\bf Proof:}
Define $W_{\phi,\Lambda\circ B}$ on $\mathbb R^N$ for 
$\phi\in  M_p^{(w)}(\mathbb Z^N)$.
\begin{eqnarray*}
W_{\phi,\Lambda\circ B}(x)& = & \sum\limits_{n\in\mathbb Z^N}
\phi(n)\Lambda\circ B(x-n)\\
& = & \sum\limits_{n\in\mathbb Z^N}\eta(n)\Lambda(Bx-n)\\
& = &W_{\eta,\Lambda}(Bx).
\end{eqnarray*}
As $\eta\in  M_p^{(w)}(\mathbb Z^N)$ (by Lemma~\ref{lem:5.4}) and 
since $\Lambda$ is a  summability kernel
 we have $W_{\eta,\Lambda}\in  M_p^{(w)}(\mathbb R^N)$. Hence 
$W_{\phi,\Lambda\circ B}\in  M_p^{(w)}(\mathbb R^N)$.
\noindent
 Similarly 
\begin{eqnarray*}
W_{\phi,\Lambda\circ B^{-1}}(x)& = &   \sum\limits_{n\in\mathbb 
Z^N}\phi(n)\Lambda(B^{-1}x-B^{-1}n)\\
& = & \sum\limits_{j=1}^q\sum\limits_{n\in B\mathbb Z^N +
p_{j}}\phi(n)\Lambda(B^{-1}x-B^{-1}n)\\
\end{eqnarray*}
where $p_1....p_{q}$ are distinct  coset representatives of 
$B\mathbb Z^N/ \mathbb Z^N\;(p_1=0).$
\begin{eqnarray*}
W_{\phi,\Lambda\circ B^{-1}}(x)& = & \sum\limits_{j=1}^q\sum\limits_{n\in \mathbb 
Z^N}\phi(Bn+p_{j})\Lambda(B^{-1}x+B^{-1}p_{j}-n)\\
& = &W_{\psi,\Lambda}(B^{-1}x)+...+W_{\psi_{p_{q-1}},\Lambda}(B^{-1}x-B^{-1}p_{q})\;\;
\end{eqnarray*}
 where $\psi_{p_i}(l)=\phi(Bl+p_j),\;i=1,2,..,q$.
 As $\psi\in  M_p^{(w)}(\mathbb Z^N)$ and $\Lambda$ is a summability 
kernel we conclude that
$W_{\phi,\Lambda\circ B^{-1}}\in  M_p^{(w)}(\mathbb R^N)$.

 Hence from Lemma~\ref{lem:6.4} and the discusssion preceeding 
 Lemma~\ref{lem:3.4} we conclude the following theorem.
\begin{theorem}
 Suppose $\Lambda\in  M_p(\mathbb R^N)$ and 
 $supp\;\Lambda\subseteq [-M,M]$; for $\phi\in  M_p^{(w)}(\mathbb Z^N)$
  define $W_{\phi, \Lambda}(\xi)= \sum
 \limits_{n\in\mathbb Z^N}\phi(n)\Lambda(\xi-n)$ on $\mathbb R^N$, 
 then $W_{\phi, \Lambda}\in  M_p^{(w)}(\mathbb R^N)$ and
$  N_p^{(w)}(W_{\phi, \Lambda})\leq C_\Lambda   N_p^{(w)}(\phi)
\|\Lambda\|_{M_p(\mathbb R^N)}$, where $C_\Lambda$ is 
a constant depending  on $\Lambda$.
\label{thm:12.4}
\end{theorem}
\section{An Application}
 As an application of Theorem~\ref{thm:12.4}, we prove a weak-type 
 version of a result proved by de Leeuw~\cite{SW}.
\begin{theorem}
 For $1<p<\infty$ , and $\epsilon>0$; let 
 $\{\phi_\epsilon\}\subseteq M_p^{(w)}(\mathbb Z)$ satisfy \\
(i) $\lim\limits_{\epsilon\rightarrow 0}\phi_\epsilon
([\frac{x}{\epsilon}])=\phi(x)$ a.e.\\
(ii) $\sup\limits_\epsilon  N_p^{(w)}(\phi_\epsilon)=K<\infty.$\\
Then $\phi\in M_p^{(w)}(\mathbb R)$ and 
$  N_p^{(w)}(\phi)\leq\sup\limits_\epsilon  N_p^{(w)}(\phi_\epsilon).$
\label{thm:13.4}
\end{theorem}
\noindent{\bf Proof:}
 For each $\epsilon>0$, define $W_{\phi_\epsilon}$ on $\mathbb R$ by
\begin{equation}
W_{\phi_\epsilon}(x)=\sum\limits_{n\in\mathbb 
Z}\phi_\epsilon(n)\chi_{[0,1)}(x-n).
\label{eqn:17.4}
\end{equation}
 As $\chi_{[0,1)}\in  M_p(\mathbb R)$ for $1<p<\infty$, from 
 Theorem~\ref{thm:12.4} we have 
 $W_{\phi_\epsilon}\in M_p^{(w)}(\mathbb R)$ and 
$  N_p^{(w)}(W_{\phi_\epsilon})\leq C  N_p^{(w)}(\phi_\epsilon)\leq CK$ 
.
 We define another function $\psi_\epsilon,$ 
for each $\epsilon>0,$  by 
$\psi_\epsilon(x)=W_{\phi_\epsilon}(\frac{x}{\epsilon})$. Then 
$\psi_\epsilon \in M_p^{(w)}(\mathbb R)$ and 
\begin{equation}
  N_p^{(w)}(\psi_\epsilon)\leq  N_p^{(w)}(W_{\phi_\epsilon})\leq CK.
\label{eqn:18.4}
\end{equation}
From (\ref{eqn:17.4}) we have
\begin{eqnarray*}
\psi_\epsilon(x)\; = \;W_{\phi_\epsilon}(\frac{x}
{\epsilon})& = & \sum\limits_{n\in\mathbb Z}\phi_\epsilon(n)
\chi_{[0,1)}(\frac{x}{\epsilon}-n)\\
& = &\phi_\epsilon([\frac{x}{\epsilon}]).
\end{eqnarray*}
So from our hypothesis  
\begin{equation}
\lim\limits_{\epsilon\rightarrow 0}\psi_\epsilon(x)=\phi(x)\;\;\; a.e.
\label{eqn:19.4}
\end{equation}
Also we have $|\psi_\epsilon(x)|<\infty\;\;\;\;$ (as 
$\sup\limits_{\epsilon,n}|\phi_\epsilon(n)|<\infty$).

 Hence from (\ref{eqn:17.4}), (\ref{eqn:18.4}) and (\ref{eqn:19.4}) 
 along with Lemma~\ref{lem:2.4} we have $\phi \in M_p^{(w)}(\mathbb R)$ 
 and $  N_p^{(w)}(\phi)\leq\lim\limits_\epsilon  
 N_p^{(w)}(\phi_\epsilon)\leq CK.$

\end{document}